\documentclass[12pt]{article}

\usepackage[margin=1in]{geometry}  
\usepackage{graphicx}           
\usepackage{amsmath}             
\usepackage{amsfonts}
\usepackage{amsthm}
\usepackage[hidelinks]{hyperref}
\usepackage{float}
\usepackage{xcolor}
\usepackage{tikz}

\usepackage{caption}
\newtheorem{theorem}{Theorem}[section]
\newtheorem{lemma}[theorem]{Lemma}

\theoremstyle{definition}
\newtheorem{definition}[theorem]{Definition}

\newtheorem{remark}[theorem]{Remark}

\newcommand{\Mcal}{\mathcal{M}}
\newcommand{\piopt}{\pi_{\mathrm{opt}}}

\newcommand{\ceilNFour}{\lceil \frac{n}{4} \rceil}
\newcommand{\floorThreeNFour}{\lfloor \frac{3n}{4} \rfloor}
\newcommand{\floorThreeNFourInline}{\lfloor 3n/4 \rfloor}

\begin{document}

\nocite{*}

\title{Transitions in optimal matchings with concave costs \\on the line}

\author{Grigorii Buklei}
\date{}

\maketitle

\begin{abstract}
We study the optimal transport problem between two finite collections of points on the real line, with the concave cost functions $|x-y|^p$ for $0 < p < 1$, which are particularly useful for economic applications. We focus on how the optimal matching changes as $p$ varies in $(0, 1)$, motivated by comparative statics in labor market job allocations. We obtain upper and lower bounds on the maximum possible number of transitions between optimal matchings. For the lower bound, we arrange points in a very specific way to link the problem to the lower envelope of line segments. For the upper bound, we use a recursive inequality derived from the geometric properties of the solution.
\end{abstract}



\section{Introduction}
\subsection{The model}

Let $x_1\leq\cdots\leq x_n$ and $y_1\leq\cdots\leq y_n$
be two collections of points on the real line, indexed in nondecreasing order.
A \textit{matching} between these collections
is a permutation $\pi$ in $S_n$, where $S_n$ denotes the set of all
permutations of $\{1, \ldots, n\}$. Given a cost function
$c:\mathbb R^2\to[0,\infty)$, the cost of the matching is
$\sum_{i=1}^n c(x_i,y_{\pi(i)})$.
In this paper, we focus on costs given by powers of the
distance, namely $c_p(x,y)=|x-y|^p$ for $p>0$. Thus, the
\textit{cost of the matching} and the \textit{minimum cost} over all
matchings are:

\vspace{-0.3em}
\[
C_p(\pi) = \sum_{i=1}^{n} |x_i-y_{\pi(i)}|^p,
\qquad F(p) = \min_{\pi \in S_n} C_p(\pi).
\]
We denote an optimal matching by $\piopt$. For $p \geq 1$, the identity permutation $\pi(i)=i$ is optimal, up to swapping coinciding points. For $p\in (0,1)$, the situation is more interesting, as the optimal matching may change.
\vspace{-1em}

\begin{figure}[h]
\includegraphics[width=1.\linewidth]{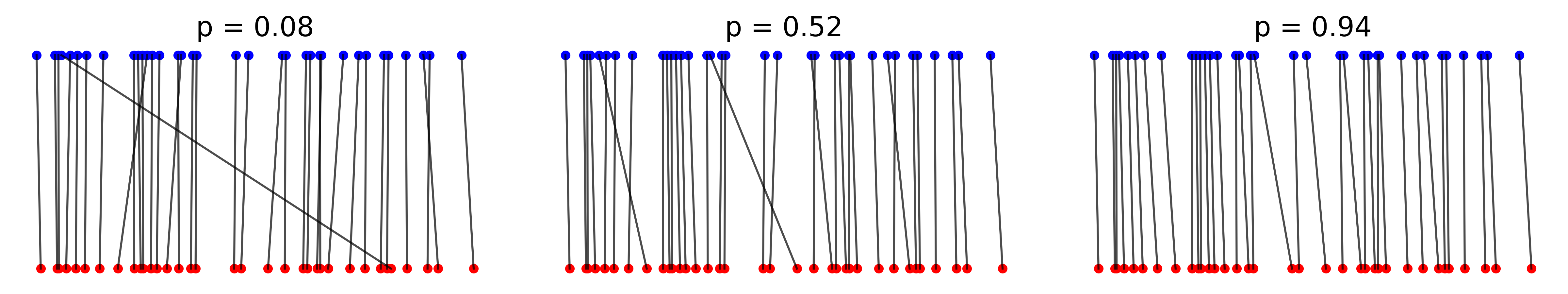}
\centering
\vspace{-2em}
\caption{Optimal matchings for different values of $p$.}
\end{figure}

\subsection{Main question}

\noindent The range $p\in(0,1)$ is especially relevant to economic applications (see Section~\ref{sec:literature-review}). In particular, it is interesting to ask how many times the optimal matching can change as $p$ varies
in $(0,1)$.

\begin{remark}\label{rmk:zero-distances}
    Pairs of coinciding points $x_i = y_j$ can be omitted in the further analysis, since they will always be matched together in the optimal matching as $|x_{i} - y_k|^p + |x_{l} - y_j|^p \geq |x_{l} - y_k|^p$ for $0 < p < 1$.
\end{remark}

\begin{remark}\label{rmk:F-analiticity-and-transitions}
For each fixed matching $\pi$, the function $C_p(\pi)$ is a finite sum of exponentials in $p$ (with zero-distance terms omitted, Remark \ref{rmk:zero-distances}), and hence is real analytic for $p>0$. Consequently, $F$ is continuous and piecewise real analytic. When optimality passes from one cost function to another, $F$ fails to be real analytic at the transition point. Conversely, every non-analytic point of $F$ corresponds to such a transition. 

Different permutations may have identical cost functions for every $p$. In this case they are optimal together whenever either is optimal, so switching between them as a chosen minimiser does not constitute a transition. Similarly, a matching that ties with the optimum only at an isolated value of $p$ does not give a transition.
\end{remark}

\vspace{0.2em}
\noindent
For a fixed configuration $({\bf x},{\bf y})$, let $V_n({\bf x},{\bf y})$ be the number of points $p\in(0,1)$ at which $F$ is not real analytic. We are interested in the largest possible number of such transitions:
\[
V_n^*=\sup_{{\bf x},{\bf y}\in\mathbb R^n}
V_n({\bf x},{\bf y}).
\]
We will show in Section \ref{sec:setup} that $V_n^*<\infty$.
Since $V_n$ takes integer values, the supremum is then attained.
\noindent
Restricting the points to $[0,1]$ would not change $V_n^*$. Indeed, any finite configuration can be translated and rescaled into $[0,1]$ without changing the optimal matchings or the number of analytic components as $p$ varies .

\subsection{Main result}

Numerical simulations with various point configurations all seemed to suggest a linear bound for $V^*_n$. In this sense, the lower bound came as a surprise.

\begin{theorem}\label{thm:main-result}
There exist constants $c_1,c_2>0$, independent of $n$, such that, for all sufficiently large $n$,

$$
c_1 \cdot n\alpha(n)\le V_n^*\le n^{c_2\log n}.
$$

\end{theorem}

\noindent Here, $\alpha(n)$ is the inverse Ackermann function \cite{ackermann1928hilbertschen}. It grows exceptionally slowly but is unbounded.
The bounds in Theorem \ref{thm:main-result} leave a substantial gap: the lower bound is superlinear, whereas the upper bound is quasi-polynomial. The true asymptotic growth of $V_n^*$ therefore remains an open question.

\begin{remark}
The upper bound comes from a recursive inequality based on the geometry of optimal matchings. For the lower bound, we arrange the points so that the matching-cost curves $p\mapsto C_p(\pi)$ reproduce the switching pattern of the lower envelope of line segments. This is where the inverse Ackermann function appears \cite{WiernikSharir1988}.
\end{remark}

\subsection*{Structure of the paper}

Section \ref{sec:literature-review} reviews the relevant literature. Section \ref{sec:setup} recalls the structural properties of optimal matchings used in the proofs. Sections \ref{sec:upper-bound} and \ref{sec:lower-bound} establish the upper and lower bounds, respectively.

\section{Literature review}\label{sec:literature-review}

Finite optimal matching is a basic instance of
optimal transport, a subject going back to Monge
\cite{Monge1781memoire} and concerned with minimising the cost of moving mass from a source to a target. It lies at the intersection of probability, analysis, and geometry \cite{Santambrogio2015, villani2009optimal}, and has attracted considerable interest for its applications in computer science, biology, and economics  \cite{galichon2016, peyré2020computationaloptimaltransport, schiebinger2019optimal}. The optimal matching problem is well understood for linear costs ($p = 1$); see Caracciolo, Erba, and Sportiello \cite{caracciolo2021OptimalMatchings}; and for convex costs, in particular $p = 2$, see Brenier \cite{brenier1991}. In the linear case, there may be multiple optimal matchings. For convex costs, we may take $\piopt(i)=i$.
\\ \\
The problem is more complicated for concave costs, which have nevertheless attracted considerable attention in recent years. For concave costs, an optimal matching does not always follow the simple ordering rule valid for convex costs. Concavity favours a long trip and a short trip over two trips of comparable length; economically, this may reflect the benefits of combining local and long-distance supply. Another application in economics is the study of worker-job matching markets; see Lindenlaub \cite{Lindenlaub2017SortingMultidimentionalTypes}, Boerma, Tsyvinski, Wang, and Zhang \cite{boerma2025compositesorting}. In this problem, the matching is an allocation of workers; it is natural to ask what happens when market conditions change, and the exponent $p$ varies. In economics, the comparison of two different outcomes, usually before and after a change in some underlying parameter, is called comparative statics; see Boerma, Ottolini, and Tsyvinski for more details \cite{boerma2026multidimensionalsortingcomparativestatics}. For example, in worker-job matching, new technology may change cost functions, and the market may adjust by reallocating workers.
\\ \\
This raises a natural question: how often do such reallocations occur when a parameter $p$ changes? This motivates our study of the number of transitions
in the optimal matching as $p$ varies between $0$ and $1$. This is challenging, since a given matching can be optimal for small $p$, then become suboptimal, and reappear as optimal later. This makes it difficult to find a useful monotone quantity as $p$ changes. McCann \cite{mccann1999exactSolutionsOnTheLine} establishes one of the main geometric properties of optimal matchings for concave costs: non-crossing, implying a layering structure (see Section \ref{sec:setup}). Relatively little is known about their geometry beyond this.
\\ \\
There has also been progress in algorithms for approximating optimal matchings. For example, Delon, Salomon, and Sobolevski \cite{DelonSalomonSobolevski2010, DelonSalomonSobolevski2012} provide methods, including a recursive algorithm, for identifying pairs that appear in an optimal matching. In \cite{CaraccioloDAchilleErbaSportiello2020}, Caracciolo, D'Achille, Erba, and Sportiello study the matching cost produced by the Dyck algorithm for different ranges of $p$. Ottolini and Steinerberger \cite{ottolini2025greedymatchingoptimaltransport} later study a greedy algorithm, establish its performance for highly concave costs ($p$ close to zero), and give bounds for both the greedy and Dyck algorithms for concave costs. Separately, Bobkov and Ledoux \cite{BobkovLedoux2019, BobkovLedoux2020} study bounds on the optimal matching cost $F(p)$ for concave cost functions. 
\\ \\
Several ingredients in the proof of Theorem \ref{thm:main-result}, including the layering structure, rely on the one-dimensional formulation and do not directly extend to higher dimensions. In dimensions $d > 1$, the geometric structure imposed by the ordering of points is lost. Nevertheless, these cases have important applications across various fields, including logistics, economics, biology, and physics \cite{Lindenlaub2017SortingMultidimentionalTypes, schiebinger2019optimal}. Some results investigate uniqueness and characterisations of optimal transport plans for concave costs \cite{GangboMcCann1996, pegon2025characterizationoptimaltransportplans}; Champion and De Pascale \cite{ChampionDePascale2011} study the existence of a transport map for general norm costs. To our knowledge, little is known about how optimal matchings vary with the cost exponent in the concave regime.
\\ \\
A probabilistic formulation -- where the points to be matched are sampled according to some probability distribution -- is useful in many applications, as it enables the analysis of average quantities rather than worst-case scenarios. Recently, the problem of matching two sets of points drawn as independent Poisson processes of equal intensity has attracted attention in the literature \cite{holroyd2011geometric, holroyd2022minimal}. These works mainly analyse structural properties of these matchings and study the existence of optimal matchings for some values of $p$; some of them treat concave costs. Various authors have studied algorithms for some concave cases and the expected matching cost \cite{CaraccioloDAchilleErbaSportiello2020, ottolini2025greedymatchingoptimaltransport}. As little is still known about the exact structure of optimal matchings for general $0 < p < 1$, a potential extension of these works is to evaluate the expected number of transitions in the optimal matching as $p$ varies, and how it evolves with $n$.

\section{Set-up and tools}\label{sec:setup}

To motivate the peculiarity of optimal matchings under concave costs, consider the following example. Let ${\bf x} = \{0.000015, 0.002700, 0.004550\}$ and ${\bf y} = \{0.002300, 0.002930, 1\}$. The figure below compares the costs of the possible matchings:

\begin{figure}[H]
\centering
\includegraphics[width=0.7\textwidth,height=.52\textheight,keepaspectratio]{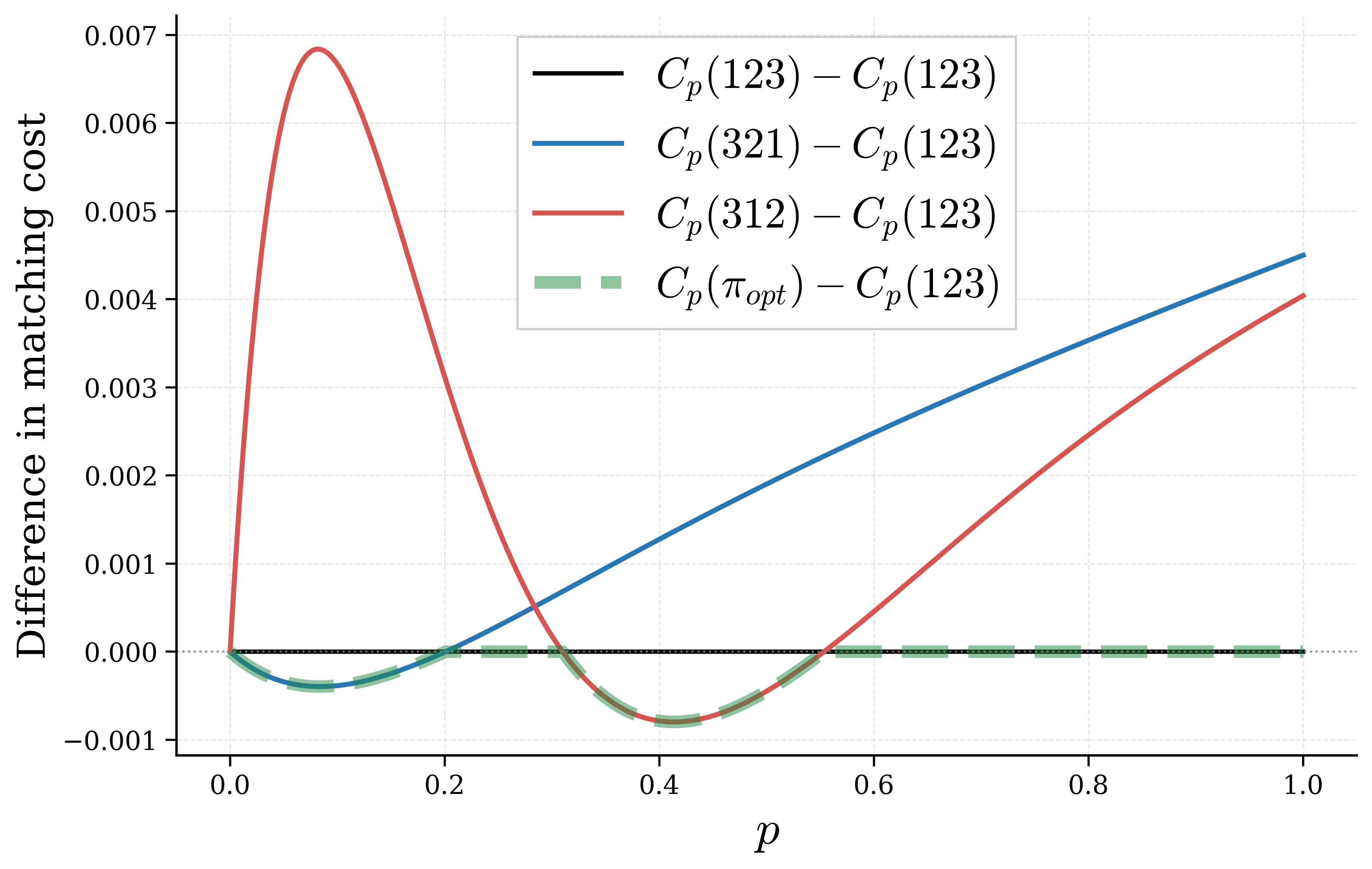}
\label{fig:counterexample-ordering}
\end{figure}

\noindent This reveals an interesting trait in the evolution of the optimal matching: $\piopt$ first equals $(3, 2, 1)$, then becomes $(1, 2, 3)$ at $p$ around $0.20$, later switches to $(3, 1, 2)$ around $p=0.31$, and becomes $(1, 2, 3)$ again after $p = 0.55$.
This rules out the hypothesis that optimal matchings can be ordered as $p$ increases.
\\ \\
Another mathematical property of concave functions is that they prioritise long and short trips over two average ones. This leads to an interesting geometrical property of the optimal solution:

\begin{theorem}[McCann, 1999 \cite{mccann1999exactSolutionsOnTheLine}]\label{thm:non-crossing}
A matching $\pi$ on the real line is non-crossing if, for every $i,j \in \{1,\ldots,n\}$, $i\neq j$, the upper semicircles whose diameters are the intervals with endpoints $x_i$, $y_{\pi(i)}$, and $x_j$, $y_{\pi(j)}$, do not intersect. For any two point sets $\mathbf{x},\mathbf{y} \subset \mathbb{R}$, the optimal matching $\piopt(p)$ is non-crossing for every $p \in (0,1).$
\end{theorem}

\begin{remark}\label{rmk:number-of-non-crossing-matching}
This greatly narrows the number of matchings that can be optimal. The number of non-crossing matchings is at most the $n$-th Catalan number $C_n=\frac{1}{n+1}\binom{2n}{n}$. 
\end{remark}
\noindent
Moreover, this restricts the set of possible matches for each point separately. Specifically, 

\begin{definition}[Layers]\label{def:layers}
Let $z_1<\dots<z_{2n}$ be all points of $(\mathbf x,\mathbf y)$ in increasing order. Put $h(0)=0$ and
$$
h(j) = 
\begin{cases}
    h(j-1)+1, \quad \text{if } z_j \text{ is an } x\text{-point} \\
    h(j-1)-1, \quad \text{if } z_j \text{ is a } y\text{-point} \\
\end{cases}$$
(so that $h(2n)=0$). Call $\ell(j)=\max\{h(j-1),h(j)\}$ the \emph{level} of $z_j$. The
\emph{layer} $L_\ell$ is the set of points of level $\ell$.
\end{definition}

\begin{figure}[H]
\centering
\includegraphics[width=1.0\textwidth,height=.52\textheight,keepaspectratio]{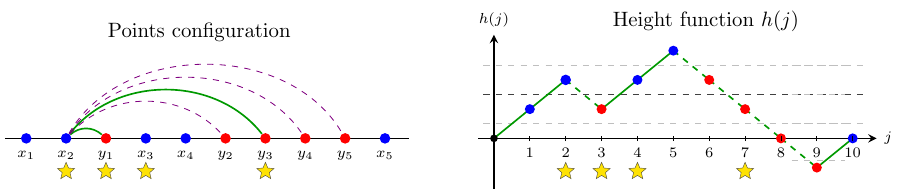}
\caption{given a point configuration (on the left) point $x_2$ can be matched only with $y_1$ or $y_3$ due to non-crossing property. This corresponds to crossing the same level set (on the right). Points within the same layer as $x_2$ are indicated with a star.}
\label{fig:layering}
\end{figure}

\begin{lemma}[Layering]\label{lem:layering}
Let $(\mathbf x,\mathbf y)$ be a configuration of $2n$ distinct points and $p\in(0,1)$.
\begin{enumerate}
  \item \emph{Arcs stay inside layers.} If $x_i$ is matched to $y_j$ in an optimal matching, then they lie in the same layer.
  \item \emph{Shape of the layers.} Each layer contains equally many $x$- and $y$-points, and each layer is an  alternating configuration of $n_\ell$ pairs, starting from $x$ or $y$.
  \item \emph{Decomposition.} Let $F_\ell(p)$ be the minimal cost of matching the points
        of $L_\ell$ among themselves. Then $F(p)=\sum_\ell F_\ell(p)$, and a matching is optimal
        for all points if and only if it pairs each point with a point of its own layer
        and is optimal inside every layer.
  \item \emph{Reduction to alternating configurations.} Let $V_n^{**}$ be the maximum of
        $V_n(\mathbf x,\mathbf y)$ over alternating configurations
        $x_1<y_1<\dots<x_n<y_n$. Then
        \[
           V_n^{**}\;\le\;V_n^{*}\;\le\;n\,V_n^{**}.
        \]
\end{enumerate}
\end{lemma}

\begin{proof}
\emph{(1).} In the optimal matching, the point $x_i$ can only be matched to a $y_j$ such that there are an equal number of $x$-type and $y$-type points between them, meaning they are in the same layer. Otherwise, interior points would be forced to match outside, breaking the non-crossing property; see Figure \ref{fig:layering}.

\emph{(2).} The function $h$ crosses the strip between heights $\ell-1$ and $\ell$ alternately upward (an $x$-point) and downward (a $y$-point). It starts at $0$, so the first crossing is upward if $\ell \ge 1$ and downward if $\ell \le 0$. Since $h(2n)=0$, up- and down-crossings are equally many.

\emph{(3).} By (1), every matched pair in an optimal matching lies within a single layer. Then, the optimal matching of all points can be obtained as unions of layerwise optimal ones. Hence, $F(p) = \sum_{\ell}F_{\ell}(p)$.

\emph{(4).} If there are $k$ layers $L_{1}, \ldots , L_{k}$ with $n_1,\ldots, n_k$ pairs in each with total number of pairs $n_1+\ldots+n_k=n$, this implies 
$$
  V^*_n \leq  \max_{\substack{k, n_1, \ldots, n_k\\ n_1+\ldots n_k=n}} \; \sum_{j=1}^k V^{**}_{n_j}\leq n V^{**}_n
$$
as both $V_n^*$ and $V_n^{**}$ are non-decreasing; see later in Lemma \ref{lem:V-0-1-increasing}. Since $V_n^{**}$ is taken over only alternating configuration of points, $V_n^{**} \leq V_n^* \leq n V^{**}_n$. In particular, this shows that the claimed upper and lower bounds in Theorem \ref{thm:main-result} for $V_n^*$ follow from that of $V_n^{**}$ up to a factor of $n$, that can be absorbed in the constant.
\end{proof}

\noindent This significantly simplifies the analysis and allows us to consider optimal matchings only within each layer, with alternating $x$ and $y$ points. To estimate the number of transitions, we compare matching costs. The difference of two matching costs is a finite sum
\vspace{-.5em}
$$f(p) = \sum_{j=1}^k a_j b_j^p, \qquad b_j>0.$$
Every exchange of optimality between two distinct cost functions occurs at a zero of their difference and the converse need not hold: the tied matchings may be suboptimal, or the equality may occur without an exchange of optimality as suggested by Remark \ref{rmk:F-analiticity-and-transitions}.
\\ \\
The following theorem bounds the number of such zeros. 

\begin{theorem}[Jameson, 2006 \cite{Jameson_2006}]\label{thm:jameson}
Let $f(p) = \sum_{j = 1}^k a_jb_j^p$ with $b_1 > b_2 > \dots > b_k > 0$, and $a_j \neq 0$. The number of real zeros $Z(f)$ of $f$ (counted with multiplicity) is bounded by the number of sign changes $S[(a_j)] $ in the sequence of coefficients $(a_1, a_2, \dots, a_k)$:
$$Z(f) \leq S[(a_j)] \leq k - 1.$$
\end{theorem}

\begin{remark}\label{rmk:number-of-changes-between-two}
Applied to the matching problem, this means that the optimiser cannot change more than $2n - 2$ times between two matchings $\pi_1$ and $\pi_2$ with different matching costs. Since the Dirichlet polynomial $C_p(\pi_1) - C_p(\pi_2) = \sum_{i=1}^n|x_i-y_{\pi_1(i)}|^p - \sum_{i=1}^n|x_i-y_{\pi_2(i)}|^p$ has one zero at $p=0$, it has at most $2n-2$ zeros in $(0, 1)$ if the two cost functions differ.
\end{remark}

\begin{remark}
Remark \ref{rmk:number-of-non-crossing-matching} notes that there are $\frac{1}{n+1}\binom{2n}{n}$ candidates for optimal matchings; combining this with Remark \ref{rmk:number-of-changes-between-two} gives an easy asymptotic upper bound on the total number of changes: $$V_n^*  = O(2^{4n}).$$
\end{remark}

\begin{remark}
This also allows us to understand the structural characterisation of optimal matching.
If all matching cost functions are different, for a generic $p$, there is a unique minimiser $\piopt(p)$, determined by which non-crossing matching has the lowest cost, and it remains so on some interval $(p_1, p_2)$. However, there are finitely many points $p^*$ where the optimiser changes. At these points, both minimisers that are optimal before and after $p^*$ are optimal. Potentially, more than two optimisers can occur at one point, but on the intervals before and after $p^*$ there is only one. \\
(picture)
\end{remark}

\noindent Another application of Dirichlet polynomials is to estimate the first values of $V_n^*$; this will be used later in the proofs.
\begin{lemma}\label{lem:V-0-1-increasing}
$V_1^* = V_1^{**} = 0, \; V_2^* = V_2^{**}= 1$, both $V_n^*$ and $V_{n}^{**}$ are non-decreasing.
\end{lemma}
\begin{proof}
    When $n=1$, the matching is unique, meaning $V_1^* = V_1^{**} = 0$.
    For two pairs of points, say, $x_1 < y_1 < x_2 < y_2$ with consecutive gaps $a, b$, and $c$, two possible matching costs are: $(a + b + c)^p + b^p$ in the nested and $a^p + c^p$ in the ordered case (see Figure \ref{fig:nested-vs-ordered}).  For a nested matching to be optimal, we need $b < \min(a, c)$. Otherwise, its cost will always be higher. \\
\vspace{-1em}

    \begin{figure}[H]
    \centering
    \includegraphics[width=0.8\textwidth,height=.52\textheight,keepaspectratio]{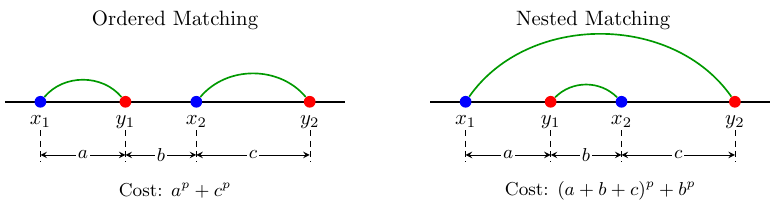}
    \caption{ordered and nested matchings in case of $n=2$.}
    \label{fig:nested-vs-ordered}
    \end{figure}

    \noindent
    The difference between their costs $(a + b + c)^p + b^p - a^p - c^p$ is then a Dirichlet polynomial with at most 2 sign changes in the coefficients: either $\underset{+}{a + b + c} \geq \underset{-}{a} \geq \underset{-}{c} \geq \underset{+}{b}$ or $\underset{+}{a + b + c} \geq \underset{-}{c} \geq \underset{-}{a} \geq \underset{+}{b}$. Since this polynomial equals zero at $p = 0$, it has no more than 1 zero in interval $p \in (0, 1)$ (Th. \ref{thm:jameson}), meaning there can be at most one change in the optimiser, $V_2^* \leq 1$.
    Taking $a = c = 0.3$, $b = 0.01$, we have $(a + b + c) \cdot b < a \cdot c$, meaning the nested will be optimal for $p$ close to 0, and $(a + b + c)+ b > a + c$, meaning the ordered will be optimal for $p$ close to 1. Hence $V_2^{**} = V_2^*\geq 1$, which proves $V_2^{**} = V_2^*=1$.

    Then we notice that given a worst-case configuration for $n$, we can place another pair of points $x_{n+1}, y_{n+1}$ sufficiently far away that it will not influence optimal matching transitions, therefore both $V_n^*$ and $V_n^{**}$ are non-decreasing.
    
\end{proof}

\section{Proof of the upper bound}\label{sec:upper-bound}

To prove a subexponential upper bound on the number of transitions, we use a recursive bound for $V_n^{**}$ and consider alternating points. Specifically, we aim to establish a bound of the form $V_n^{**} \leq C n^{\beta} \cdot V_{\floorThreeNFourInline}^{**}$. By Lemma \ref{lem:V-0-1-increasing}.(4), $V_n^*$ then satisfies the same bound up to a change in the constant, which later implies the upper bound. Although the constant $C$ and $\beta$ could be improved, doing so would not change the asymptotic class of the bound. We therefore use some deliberately loose estimates for simplicity.
\begin{remark}
In this section, when referring to arc length, we mean the number of $x$-type points (or, equivalently, $y$-type points) it embraces, since relative distances between points will not be used in the proofs. We will also operate with combined set of ordered points retaining their labels, and using a shared indexing:
$$0 \; < \; z_1 = x_1 \; < \; z_2 = y_1  \; < \; \cdots \; < \; z_{2n-1}=x_n \; < \; z_{2n} = y_n \; < \; 1. $$
\end{remark}

\noindent
To establish the recursive bound, we introduce some terms:
\begin{definition}
    A \textit{block of size $r$} is a set $I = \{k, k+1, \ldots, k+2r -1\} \subset \{1, \ldots, 2n\}$ of consecutive elements. This block has an even length; hence, there are the same number of $x$ points (with odd indices) and $y$ points (with even indices) in $\{z_i\}_{i \in I}$. Then also define:
    \begin{itemize}
        \item $F^A$, minimal cost of matchings of points with indices in $A$, where $A \subseteq \{1, \ldots, 2n\}$ a set of combined point indices that is balanced, having same number of $x$ and $y$ points.
        $$F^A(p)= \min_{\substack{\sigma}} \sum_{\substack{i\in A\\ i\text{ is odd}}}|z_i-z_{\sigma(i)}|^p, \quad \text{where } \sigma \text{ is bijective, } \sigma(i) \in A.$$
        \item $\Gamma$ -- set of all blocks $I = \{k, k+1, \ldots, k+2r -1\}$ of size $r \in \big[\ceilNFour; \floorThreeNFour\big]$. \\[0.5em]
        We write $I^c $ for $\{1, \dots, 2n\} \setminus I$. \
        \textit{Note:} $|\Gamma| \leq \sum_{r=1}^{n}(2n-2r+1) =n^2$
\end{itemize}
\end{definition}

\begin{lemma}\label{lem:fraction-consecutive-points}
For every non-crossing matching $\pi$ between $n$ pairs of points $x$ and $y$, $n \geq 2$, there exists a block 
$$I = \{k, k+1, \dots, k + 2r-1\},\quad \text{of size } 
\left\lceil\frac{n}{4}\right\rceil
\le r\le
\left\lfloor\frac{3n}{4}\right\rfloor$$
such that all points in $I$ are matched among themselves. 
\end{lemma}
\begin{proof}
    Since $\pi$ is non-crossing, any two arcs are nested or disjoint. Call arcs that lie inside no other arc \textit{outermost}. In a non-crossing matching of all points the longest arc is outermost; say it has length $l_0$ (meaning it has $l_0 - 1$ points of each type inside). Consider three cases:
    \begin{enumerate}
        \item if $l_0 < \lceil\frac{n}{4}\rceil$. Then all outermost arcs have lengths $h_1, \dots, h_p$ smaller than $\ceilNFour$, and altogether cover the whole line. Select first $m$ of them such that $\sum_{i=1}^m h_i \geq \ceilNFour$. Then, 
        $$\sum_{i=1}^m h_i = \sum_{i=1}^{m-1} h_i + h_m \leq \Big\lceil\frac{n}{4}\Big\rceil - 1 +  \Big\lceil\frac{n}{4}\Big\rceil - 1 = 2 \Big\lceil\frac{n}{4}\Big\rceil - 2 \leq \frac{n-1}{2} \leq \Big\lfloor\frac{3n}{4}\Big\rfloor$$ 
        \item if $l_0 \in \big[\ceilNFour; \floorThreeNFour \big]$. In this case, the interval set $I$ of all points that the arc of length $l_0$ embraces is self-matched and has the required length.
        \item if $l_0 > \floorThreeNFour$. Then, consider all the points that are inside that arc. Let the longest arc be of length $l_1$. We will continue the procedure until, at some stage $s$, $l_{s+1} \leq \floorThreeNFour < l_{s}$.
        Then, if $l_{s+1} \in \big[ \ceilNFour; \floorThreeNFour \big]$, the proof is finished as in case (2). Otherwise, if $l_{s+1} < \ceilNFour$, then arc $l_s$ contains at least $\floorThreeNFour$ pairs of points, and the longest arc is shorter than $\ceilNFour$, so the argument concludes as in case (1).
    \end{enumerate}
\end{proof}

\vspace{-1.3em}
\begin{figure}[H]
\centering
\includegraphics[width=0.9\textwidth,height=.52\textheight,keepaspectratio]{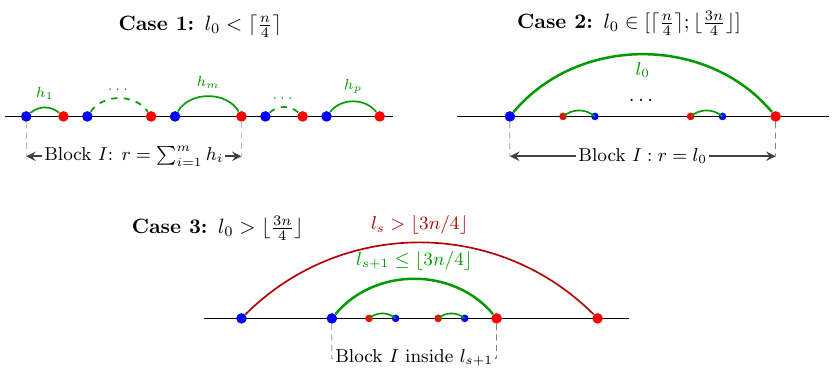}
\end{figure}

\noindent
In these notations, using Lemma \ref{lem:fraction-consecutive-points}, 
\begin{equation}\label{eq:F_n_min_definition}
F(p) = \min_{I \in \Gamma} \big( F^{I}(p) + F^{I^c}(p) \big).
\end{equation}

\begin{lemma}\label{lem:fixed_interval_V_n}
For $I \in \Gamma$, $F^{I}(p)+F^{I^c}(p)$ has at most $2 \cdot V_{\lfloor3n/4\rfloor}^{**}$ non-analytic points, and the function is analytic between them.
\end{lemma}
\begin{proof}
Given a block $I$, the functions $F^{I}(p)$ and $F^{I^c}(p)$ are optimal matching cost functions of $r \in \big[\ceilNFour; \floorThreeNFour\big]$ points. Because $V_n^{**}$ is increasing (see Lemma \ref{lem:V-0-1-increasing}), each function has at most $V_{\lfloor 3n/4 \rfloor}^{**}$ non-analytic points. Between these points, $F^{I}(p)$ and $F^{I^c}(p)$ are analytic, so their sum is analytic as well. Therefore, $F^{I}(p)+F^{I^c}(p)$ has at most $2 \cdot V_{\lfloor 3n/4\rfloor}^{**}$ non-analytic points.
\end{proof}

\begin{lemma}\label{lem:local_bounds_V_n}
Denote by $w_1 < \dots < w_k \in (0, 1)$ all non-analytic points of the functions $F^{I}(p)+F^{I^c}(p)$ for all $I \in \Gamma$. Then
\begin{enumerate}
\item $k \leq 2n^2 \cdot  V_{\lfloor3n/4\rfloor}^{**}$
\item Add $w_0 = 0$ and $w_{k+1} = 1$. On each interval $(w_s, w_{s+1})$, the optimal matching cost function $F(p)$ has fewer than $n^5$ non-analytic points
\end{enumerate}
\end{lemma}
\begin{proof}
\textit{(1.)}
$\Gamma$ consists of at most $n^2$ blocks. For each block $I$, $F^{I}(p)+F^{I^c}(p)$ has no more than $2 \cdot V_{\lfloor3n/4\rfloor}^{**}$ non-analytic points, and the function is analytic between them (Lemma \ref{lem:fixed_interval_V_n}). This proves the first statement.

\textit{(2.)} Optimal matching cost function $F(p)$ is the minimum over all $I \in \Gamma$. For each pair of different $I_1, I_2 \in \Gamma$, the minimum 
$$\min_{j \in \{1, 2\}} \bigl(F^{I_j}(p)+F^{I_j^c}(p)\bigr),$$
\noindent
has at most $2n-2$ non-analytic points on $(w_s, w_{s+1})$; see Remark \ref{rmk:number-of-changes-between-two}. Since there are at most $\binom{|\Gamma|}{2} \le n^2 \cdot (n^2-1) / 2$ ways to choose such a pair, $F(p)$ cannot have more than: $2n \cdot \binom{|\Gamma|}{2} < n^5$ non-analytic points in a fixed $(w_s, w_{s+1})$.
\end{proof}

\begin{lemma}\label{lem:recursive-bound-Vn}
For $n \geq 3$, $V_n^{**} \leq 3 n^{7} \cdot V_{\floorThreeNFourInline}^{**} $.
\end{lemma}
\begin{proof}
From the definition (\ref{eq:F_n_min_definition}) of $F(p)$:
$$F(p) = \min_{I \in \Gamma} \big( F^{I}(p) + F^{I^c}(p) \big),$$
its non-analyticity  can happen in 2 ways: \textit{(1)} another block $I \in \Gamma$ becomes  optimal or \textit{(2)} for the same minimiser $I$ an optimal matching changes in $I$ or $I^c$.

The number of changes of \textit{(2)} corresponds to the first part of Lemma \ref{lem:local_bounds_V_n}, and there are at most $2 n^2 \cdot V_{\lfloor3n/4\rfloor}^{**}$ such points. These points create at most $2n^2 \cdot V_{\lfloor3n/4\rfloor}^{**} +1$ intervals, and each interval contains at most $n^5$ points, by the second part of Lemma \ref{lem:local_bounds_V_n}. Therefore, the maximum number of non-analytic points of $F(p)$ with alternating points:
\begin{align*}
    V_n^{**} &< 2n^2 \cdot V_{\floorThreeNFourInline}^{**} + \big(2n^2 \cdot V_{\floorThreeNFourInline}^{**} + 1\big)\cdot n^5  < 3n^{7} \cdot  V_{\floorThreeNFourInline}^{**}
\end{align*}
\end{proof}

\begin{lemma}\label{lem:upper-bound}
For $n \geq 2$,  $V_n^{*} \leq n^{c_2 \log n}$, where $c_2$ is independent of $n$.
\end{lemma}
\begin{proof}
Define $n_1 := n$, and $n_{k+1} = \lfloor \frac{3}{4} n_k \rfloor$ until $\lfloor \frac{3}{4} n_K \rfloor = 2$. Then,
\begin{equation}
V_n^{**} \leq 3n_1^{7} \cdot 3n_2^{7} \cdot \dots \cdot 3n_K^{7} \cdot  V_2^{**} = \prod_{k = 1}^{K} 3n_k^{7} \leq (3n^{7})^K
\end{equation}
since $V_2^{**}=1$ from Lemma \ref{lem:V-0-1-increasing}.
\noindent
Since $n_k \leq n(3/4)^{k-1}$, we have $K=O(\log n)$.
Consequently, using Lemma \ref{lem:layering},
$$V_n^* \leq nV_n^{**} \le ne^{C(\ln n)^2} =n^{c_2\ln n}$$
for a constant $c_2>0$ independent of $n$.
\end{proof}

\noindent This proves the upper bound from Theorem \ref{thm:main-result}.

\section{Proof of the lower bound} 
\label{sec:lower-bound}

The proof of the lower bound links the combinatorial construction of lower envelope of segments with matching construction. The structure of the proof (informally) goes as follows:
\\ \\
\emph{Step 1. Segments construction.} \\Construct $m$ line segments $s_i$ whose lower envelope changes $k=\Theta\bigl(m \alpha(m)\bigr)$ times.
\vspace{-0.5em}
$$\min_{i = 1}^m s_i(p)$$

\noindent
\emph{Step 2. Segment approximation with exponentials.}
\\Approximate each segment by a sum of four exponentials
\vspace{-0.5em}
$$h_i(p) = \sum_{j=1}^{4} f_{ij}e^{g_jp}$$

\noindent
\emph{Step 3. Matching construction.}
\\ Arrange points on the line such that matching-cost differences reproduce these exponential functions, up to small error.
\vspace{-1em}
$$\min_{i= 1}^n C_p(\Mcal_i)$$
\\ \\
\noindent
We describe an informal overview of how the following lemmas fit together. Take a value $p_r$ from Lemma $\ref{lem:analytic-points-in-area}$, and suppose the unique segment attaining the minimum at this value is $s_{i_r}$. Thus,
$$\min_{i} s_i(p_r) = s_{i_r}(p_r).$$
Using Lemma \ref{lem:exponential-approximation}, we approximate segments with sums of exponential functions $h_i(p)$ in a way that the approximation preserves the strict separation between the minimiser and the other functions:
$$\min_i h_i(p_r)=h_{i_r}(p_r).$$
Later, in Lemmas \ref{lem:selection-alpha-L}, \ref{lem:choose-r-q}, \ref{lem:matchings-cost-difference}, we allocate points $\bf x, y$ in a way that, up to a normalisation, the matching costs become functions $h_i$ to a sufficiently small error:
$$\frac{C_p(M_{4i-3})-C_p(M_1)}{2\alpha p}+h_1(p)=h_i(p)+\gamma_i(p),\qquad |\gamma_i(p)|<\delta/2.$$
Then, in Lemma \ref{lem:points-as-h-i} we match the number of points $p_r$, with subsequently different minimisers, to number of minimisers of optimal matchings:
$$\operatorname*{argmin}_{i} C_{p_r}(\Mcal_{4i-3}) = \operatorname*{argmin}_i h_i(p_r),$$
and prove the lower bound in Lemma \ref{lem:lower-bound}.

\subsection{Constructing segments with many transitions}

In this section, we use the result of Wiernik and Sharir \cite{WiernikSharir1988} to construct $m$ segments whose lower envelope has $\Theta(m \alpha(m))$ transitions. 
\begin{lemma}\label{lem:analytic-points-in-area}
For any $m \in \mathbb{N}$, there exist $m$ segments $s_1, \ldots, s_m$ in $[1/3, 2/3] \times [0,1]$, points $\frac{1}{3}\leq p_1 < p_2 <\cdots < p_{k}\leq \frac{2}{3}$, with $k \geq cm\alpha(m)$, and $\delta>0$ such that, for each $r$, the interval $(p_r-\delta, p_r+\delta)$ contains a unique minimising segment $s_{i(r)}$, and these minimisers change between consecutive $p_r$. Moreover, any other segment defined on the interval is at least $2\delta$ larger.
\end{lemma}
\begin{proof}
In \cite{WiernikSharir1988}, Wiernik and Sharir show that for any number $m \in \mathbb{N}$ we can select $m$ segments $s_i$ of the lines $a_i+b_ip$ on $[c_i, d_i]$ such that their lower envelope: $Y_S(p) = \min_{i} s_i$ has $\Theta(m\alpha(m))$ pieces.

Without loss of generality, we may assume that all segments are contained in $[1/3, 2/3] \times [0,1]$, since any configuration can be rescaled and mapped to this region, which preserves the set of minimisers. This lets us exclude the case where $p$ is close to $0$ or $1$ in the further arguments.

Consider, for each analytic component of the lower envelope $Y_S$, a midpoint, and denote these points by $p_1 < \cdots < p_k$. By construction, each $p_r$ lies in $[1/3, 2/3]$. Furthermore, the number of such points satisfies $k \geq c m \alpha(m)$ for some constant $c > 0$ independent of $m$.

It remains to show that these points can be separated. Since each analytic component is an interval or segment, for each $p_r$ there exists $\delta_1 > 0$ such that $(p_r - \delta_1, p_r + \delta_1)$ is contained within the same analytic component. Moreover, since the minimiser is unique at each $p_r$ and the number of segments is finite, there exists $\delta_2 > 0$ such that, throughout $(p_r - \delta_1, p_r + \delta_1)$, any other segment remains at least $\delta_2$ above the minimising segment. We also note that, since there are finitely many segments and sample points, there exists $\rho>0$ such that
$$\operatorname{dist}(p_r,[c_i,d_i]) \geq \rho \qquad\text{whenever }p_r\notin[c_i,d_i]. 
$$
Setting $\delta = \min\{\delta_1,\delta_2/2,\rho\}$ completes the proof of the lemma.
\end{proof}

\subsection{Approximation by exponentials}

\begin{lemma}\label{lem:exponential-approximation}
Given $m$ segments, for any $\delta>0$ we can approximate them by functions $\{h_i(p) = \sum_{j=1}^4 f_{ij}e^{g_jp}\}_{i=1}^m$ in a way that, on each segment, they stay within $\delta/2$ of the segment, and at distance $\delta/2$ from the segment's domain they are greater than 2.
\end{lemma}
\begin{proof}
Let the segments be defined by coefficients $\{a_i, b_i, c_i, d_i\}_{i}^m$, so that segment $s_i$ is part of the line $a_i + b_ip$ on the interval $[c_i, d_i]$. Then define
$$A = \max_{i=1}^m \big( |a_i| + |b_i| \big), \quad B = \max_{i=1}^m \big( |b_i| \big).
$$
Select $g_1 = 0, g_2 = \nu, g_3 = K, g_4 = -K$ such that
$$\nu = \min \Big\{1, \frac{\delta}{2eB}\Big\},
\quad K = \frac{6}{\delta} \ln{\Big(\max \Bigl\{\frac{8}{\delta}, 4 + A + \frac{\delta}{4}\Bigr\}\Big)}.
$$
Set
$$f_{i1} = a_i - \frac{b_i}{\nu}, \quad f_{i2} =  \frac{b_i}{\nu}, \quad f_{i3} = e^{-K(d_i + \delta/6)}, \quad f_{i4} =e^{K(c_i - \delta/6)},$$
and
$$h_i(p) =
\underbrace{a_i + \frac{b_i}{\nu}(e^{\nu p} -1)}_{L_i(p)} +
\underbrace{e^{K(p - d_i - \delta/6)}}_{V_i(p)} +
\underbrace{e^{-K(p - c_i + \delta/6)}}_{W_i(p)}
.$$
Since $|e^x - 1 - x| \leq \frac{x^2}{2}e^{|x|}$ (Taylor expansion), then
\begin{equation}\label{eq:L_i(p)-estimate}
    |L_i(p) - (a_i + b_ip)| =
    \frac{|b_i|}{\nu}|e^{\nu p} - 1 - \nu p| \leq
    \frac{|b_i|}{\nu}\frac{\nu^2p^2}{2}e^{\nu p} \leq
    \frac{|b_i| \nu e^\nu}{2} \leq \frac{\delta}{4}.
\end{equation}
Moreover,
\begin{itemize}
    \item Inside $[c_i, d_i]$.
    By the choice of $K$, $e^{K \delta/6} \geq \frac{8}{\delta}$, $e^{K \delta/6} \geq 4 + A + \frac{\delta}{4}$. Then, in $[c_i, d_i]$:
    \[
    \begin{cases}
        V_i(p) \leq e^{-K \delta/6}\\
        W_i(p) \leq e^{-K \delta/6}
    \end{cases} \quad
    \Rightarrow \quad
    \|h_i(p) - (a_i + b_ip)\|_{\infty}
    \;\; \leq \;\;
    \frac{\delta}{4}+\frac{\delta}{8}+\frac{\delta}{8} = \frac{\delta}{2}.
    \]
    \item Outside $[c_i, d_i]$. Consider $p \geq d_i + \delta/2$. Then,
    \[
    \begin{cases}
        V_i(p)\geq e^{K \delta/6} \\
        W_i(p) \geq 0 \\
        L_i(p) \geq -A - \frac{\delta}{4} \;\; (\text{using eq. } (\ref{eq:L_i(p)-estimate}))
    \end{cases}
    \Rightarrow h_i(p) \geq A + \frac{\delta}{4} + 4 - A - \frac{\delta}{4} > 2
    \]
    The case $p \leq c_i - \delta/2$ is symmetric, with $W_i$ in place of $V_i$.
\end{itemize}
\end{proof}

\subsection{Constructing points on the line}

\subsubsection{Placing the points}\label{subsec:placing-points}

We construct two sets of ordered points $x_1 < y_1 < \ldots < x_n < y_n$, with $\{x_i, y_i\}_{i=1}^{n} \subset [0, 1]$. The construction we give may exceed 1 on the right, but since the number of points is finite, we can rescale the whole configuration to \([0, 1]\), preserving the ordering and the count of points.

The construction relies on placing $x$-points at distances $2r_i$, and placing $y$-points roughly in between, controlled by parameters $r_i$ and $\alpha$, except for the last one, which is far away. The distances between points then will be:
\begin{align*}
    &y_i -x_i = r_i (1+\alpha q_i) &\text{ for }& i\in\{1, \dots, n-1\}, \\
    &x_{i+1} -y_i = r_i (1-\alpha q_i) &\text{ for }& i\in\{1, \dots, n-1\}, \\
    &y_n -x_n = L,  && 
\end{align*}
where we require $\alpha$ to be small enough to keep the distances positive. For arbitrary choices of $r_i > 0$ and $q_i \in \mathbb R$, one can arrange the points so that the optimal matching is determined by pairing a single element $y_n$. Thus, the matching is simple to describe, while there is still substantial flexibility.

\begin{definition}
    We will denote the matching where the last point $y_n$ is matched to the point $x_i$, and others are matched in order, as $\mathcal{M}_i$. So, $\mathcal{M}_i(j) = j$ for $j < i$, $\mathcal{M}_i(j) = n$ for $j = i$, and $\mathcal{M}_i(j) = j-1$ for $j > i$, see on the Figure \ref{fig:points-allocation}.
    \begin{figure}[H]
    \centering
    \includegraphics[width=0.95\textwidth,height=.52\textheight,keepaspectratio]{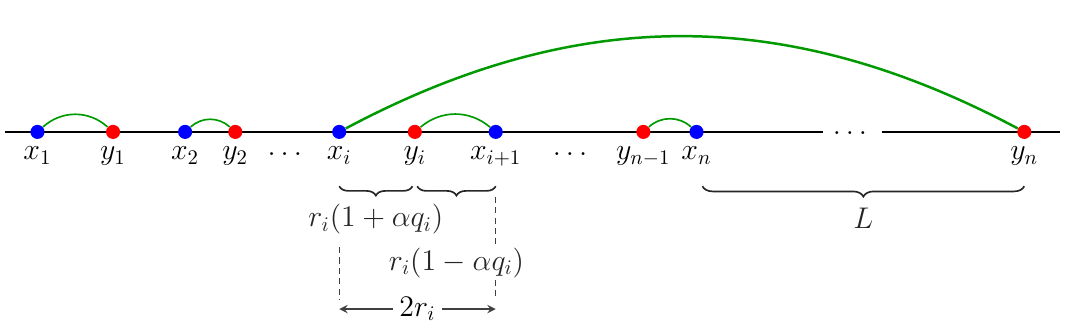}
    \caption{Matching $\Mcal_i$.}
    \label{fig:points-allocation}
    \end{figure}
\end{definition}

\subsubsection{Restricting the candidate matchings}
In this section we select such $\alpha$ and $L$ so that the choice of optimal matchings will be limited to some $\Mcal_i$, and prove a technical inequality which will be used later.

\begin{lemma}\label{lem:selection-alpha-L}
Fix $r_j > 0$ and $q_j \in \mathbb R$, $j = 1, \ldots, n - 1$, and constants $c, \delta > 0$. There exist $\alpha > 0$ and $L \geq 1$, with $\alpha \max_j |q_j| \leq 1/2$, such that:
\begin{enumerate}
\item For every $p\in[1/3,2/3]$ and every $i=1,\ldots,n$, the matching $\mathcal M_i$ is the unique minimum among matchings in which $y_n$ is matched to $x_i$.
\item
\[
c\alpha+
\sup_{\substack{p\in[1/3,2/3] \\ 1 \leq i \leq n}}
\frac{\bigl||x_i - y_n|^p-|x_1 - y_n|^p\bigr|}
     {2p\alpha}
<\frac{\delta}{2}.
\]
\end{enumerate}
\end{lemma}
\begin{proof}
We will establish these conditions one by one, and at each step we can make $\alpha$ smaller or $L$ larger in a way that does not affect the previous inequalities.

\emph{(1).} 
    We will prove that if $y_n$ is matched to $x_i$, then $\Mcal_i$ is the optimal matching. By the non-crossing condition $x_1,\ldots, x_{i-1}$ are matched to $y_1,\ldots, y_{i-1}$ and  $x_{i+1},\ldots, x_{n}$ are matched to $y_i,\ldots, y_{n-1}$. We will prove that the ordered matching is the unique optimal matching for the first set of points; the argument for the second set is analogous.
    
Suppose that $y_k$ is matched to some $x_j$ with $j,k<i$, then its cost is at least $|r_k(1 \pm \alpha q_k)|^p$. If $j \neq k$, an additional term arises, so the cost is at least $|r_k(1 \pm \alpha q_k) + r_m(1 \pm \alpha q_k)|^p$. Denote by $\pi$ the matching between these points, and $\pi_{ord}$ the ordered matching. The difference in cost is then given by:
\begin{align*}
    C_p(\pi) -& C_p(\pi_{ord}) \geq \\ 
    &|r_k(1 \pm \alpha q_k) + r_m(1 \pm \alpha q_k)|^p + \sum_{t = 1, t \neq k}^{i-1} |r_t(1 \pm \alpha q_t)|^p -
    \sum_{t = 1}^{i-1} |r_t(1 \pm \alpha q_t)|^p  = \\
    &|r_k(1 \pm \alpha q_k) + r_m(1 \pm \alpha q_k)|^p - |r_k(1 \pm \alpha q_k)|^p = \\[0.7em]
    &|r_k + r_m|^p - |r_k|^p + O(\alpha). 
\end{align*}
\noindent
Therefore, for fixed $r_j > 0$ and $q_j$, it is possible to choose $\alpha$ sufficiently small so that any matching $\pi$ has strictly greater cost than $\pi_{ord}$. Since the number of possible matchings is finite, we may select the minimal admissible value of $\alpha$, subject to the constraint $\alpha \max_i |q_i| < 1$, as required by the points construction.

\emph{(2).} Given any $c > 0$, select $\alpha > 0$ sufficiently small so that $c\alpha < \delta/4$, while maintaining the conditions established in part (1). By definition,
    $$ |x_1-y_n|^p-|x_i-y_n|^p =
    \left(L+2\sum_{j=1}^{n-1}r_j\right)^p - 
    \left(L+2\sum_{j=i}^{n-1}r_j\right)^p.
    $$
    For fixed $r_j$ and $\alpha$, since $p\in[1/3,2/3]$, the mean value theorem implies that
    $$ \sup_{p,i} \frac{\big||x_1-y_n|^p-|x_i-y_n|^p\big|}{2p\alpha} \to 0 \quad \text{as }L\to\infty.
    $$
    Thus, by taking $L$ sufficiently large, we obtain
    $$
    c\alpha+ \left| \frac{|x_1-y_n|^p-|x_i-y_n|^p}{2p\alpha} \right| <\frac{\delta}{4}+\frac{\delta}{4} =\frac{\delta}{2}.$$
\end{proof}

\subsubsection{Matching cost differences}
Here we establish the relation between matching cost, normalised matching cost differences, and exponential functions.

\begin{lemma}\label{lem:choose-r-q}
    Let
    $$ h_i(p)=\sum_{j=1}^4 f_{ij}e^{g_jp},
    \qquad i=1,\ldots,m, $$
    and set $n=4m-3$. There exist $r_j>0$ and $q_j\in\mathbb R$, $j=1,\ldots,n-1$, such that, with
    
    $$P_\ell(p):=h_1(p)+\sum_{j=1}^{\ell-1}r_j^p q_j,
    \qquad 1\leq\ell\leq n, $$
    we have
    $$ P_{4i-3}(p)=h_i(p), \qquad i=1, \ldots, m, $$
    and
    $$P_{4i-3+t}(p)\geq\min\{h_i(p),h_{i+1}(p)\},
    \qquad i=1,\ldots,m-1,\quad t=1,2,3. $$
    The latter inequality is strict whenever $h_i(p)\neq h_{i+1}(p)$. In particular,
    $$\min_{1\leq\ell\leq n}P_\ell(p) =
    \min_{1\leq i\leq m}h_i(p).$$
\end{lemma}

\begin{proof}
    For $i=1,\ldots,m-1$, put
    $$ s_{ij}=f_{i+1,j}-f_{ij}. $$
    Choose a permutation $\sigma_i$ so that the positive $s_{ij}$ come first, followed by the zero coefficients and then the negative coefficients. For $t=1,\ldots,4$, define
    $$ q_{4(i-1)+t}=s_{i,\sigma_i(t)}, \qquad
    r_{4(i-1)+t}=e^{g_{\sigma_i(t)}}. $$
    Then
    $$ \sum_{t=1}^4 r_{4(i-1)+t}^p q_{4(i-1)+t} 
    = h_{i+1}(p) - h_i(p). $$
    Since $P_1=h_1$, we obtain $P_{4i-3}=h_i$. For each fixed $p$, the four increments have the same signs as the corresponding coefficients $s_{ij}$. The partial sums therefore first increase through the positive increments, remain constant through the zero increments, and then decrease through the negative increments. Therefore, every intermediate partial sum is strictly larger than the smallest of the endpoints.
\end{proof}

\begin{lemma}\label{lem:matchings-cost-difference}
    Consider the configuration $x_1<y_1<\dots<x_n<y_n$ of Section \ref{subsec:placing-points}, satisfying Lemmas \ref{lem:selection-alpha-L}, \ref{lem:choose-r-q}. Then,
    $$C_p(\mathcal{M}_i) - C_p(\mathcal{M}_1) = 
    2\alpha p \cdot \sum_{j=1}^{i-1}r_j^pq_j + \big(|x_i - y_n|^p - |x_1 - y_n|^p\big)+ O(\alpha^2)$$
\end{lemma}

\begin{proof}
Under the matchings $\Mcal_i$ and $\Mcal_1$, points $x_{i+1}, \ldots, x_n$ are matched to $y_i, \ldots, y_{n-1}$ respectively, in the same way.
Then the cost difference is:
    \begin{align*}
          &|x_1 - y_1|^p + \ldots  + |x_{i-1} - y_{i-1}|^p \quad+|x_i - y_n|^p \\
        - &|y_1 - x_2|^p - \ldots  - |y_{i-1} - x_i|^p \qquad-|x_1 - y_n|^p \\[0.5em]
        = &(r_1(1 + \alpha q_1))^p + \ldots + (r_{i-1}(1 + \alpha q_{i-1}))^p \quad+|x_i - y_n|^p \\
        - &(r_1(1 - \alpha q_1))^p - \ldots - (r_{i-1}(1 - \alpha q_{i-1}))^p \quad-|x_1 - y_n|^p \\[0.5em]
        = &\sum_{j=1}^{i-1} r_j^p\Big((1 + \alpha q_j)^p - (1 - \alpha q_j)^p\Big)  + |x_i - y_n|^p - |x_1 - y_n|^p \\
        &= 2\alpha p \cdot \sum_{j=1}^{i-1}r_j^pq_j + O(\alpha^2) + |x_i - y_n|^p - |x_1 - y_n|^p.
    \end{align*}
\end{proof}

\begin{lemma}\label{lem:points-as-h-i}
Given functions 
\vspace{-0.7em}
$$h_i(p) = \sum_{j=1}^4 f_{ij}e^{g_jp}, \quad i \in \{1, \ldots m\},$$ 
suppose that:
    \begin{enumerate}
        \item There exist points $p_1<\cdots<p_k$ in $[1/3,2/3]$ such that, for each $r\in\{1,\ldots,k\}$, the function values $h_i(p_r)$ have a unique minimiser $i_r$. These minimisers differ at consecutive points: $i_r\ne i_{r+1}$ for $r=1,\ldots,k-1$.
        \item There exists $\delta > 0$ such that
        $$h_i(p_r) \geq h_{i_r}(p_r)+\delta$$
        for every $r$ and every $i \neq i_r$.
    \end{enumerate}
\noindent
Then there is a configuration of $n=4m-3$ pairs of points such that the optimal matching changes at least $k-1$ times for $p\in[1/3, 2/3]$.
\end{lemma}

\begin{proof}
    First, we select $n = 4m-3$ pairs of points as in Lemma \ref{lem:choose-r-q}. By that lemma and finiteness, we may decrease $\delta$ so that $P_\ell(p_r)\ge h_{i_r}(p_r)+\delta$ for every $r$ and $\ell \ne 4i_r - 3$. Then, note that adding a smooth function or multiplying each term by a positive smooth function does not change the argmin. Hence,
    \begin{align*}
        \operatorname*{argmin}_{i=1}^{n} \; C_p(\Mcal_i)
        \quad \Leftrightarrow \quad
        \operatorname*{argmin}_{i = 1}^{n} \; \frac{C_p(\Mcal_i)-C_p(\Mcal_1)}{2p\alpha} + e^{g_1p}f_{11} + e^{g_2p}f_{12} + e^{g_3p}f_{13} + e^{g_4p}f_{14}\
    \end{align*}
    \noindent
    Moreover, from Lemma \ref{lem:matchings-cost-difference} we know that for small $\alpha$, 
    $$\frac{C_p(\Mcal_i)-C_p(\Mcal_1)}{2p\alpha} = \sum_{j=1}^{i-1} r_j^p q_j + \underbrace{O(\alpha) + \frac{|x_i - y_n|^p - |x_1 - y_n|^p}{2p\alpha}}_{\gamma_i},$$ 
    and, from Lemma \ref{lem:selection-alpha-L}, that $|\gamma_i|$ can be limited by $\delta/2$:
    $$\frac{C_p(\Mcal_i)-C_p(\Mcal_1)}{2p\alpha} = \sum_{j=1}^{i-1} r_j^p q_j + \gamma_i, \quad |\gamma_i| < \delta/2.$$ 
    \noindent
    From Lemma \ref{lem:choose-r-q} we also know that taking these $\bf r, q$, we get 
    \begin{align*}
        &e^{g_1p}f_{11} + e^{g_2p}f_{12} + e^{g_3p}f_{13} + e^{g_4p}f_{14} + \sum_{j=1}^{4i-4}r_j^p q_j \\ &= 
        e^{g_1p}f_{i1} + e^{g_2p}f_{i2} + e^{g_3p}f_{i3} + e^{g_4p}f_{i4} =  h_i(p).
    \end{align*}
    Then, at points $p_1, \ldots, p_k$ matching $\Mcal_{4i-3}$ is optimal if and only if $i = \operatorname*{argmin}_i h_i(p)$ since we know that at these points
    $$\operatorname*{argmin}_{i=1}^{n=4m-3}\; C_p(\Mcal_i) 
    \;\; = \;\;
    \operatorname*{argmin}_{i=1}^{m}\; C_p(\Mcal_{4i-3}) 
    \;\; = \;\; 
    \operatorname*{argmin}_{i=1}^m\; h_i(p) + \gamma_i.$$
    By the lemma’s assumption, at each point $p_r$ the minimum value is separated from every other value by at least $\delta$. Since there are finitely many indices, we can choose $\alpha$ small enough that $|\gamma_i|<\delta/2$ for every $i$. Thus, the unique minimiser at each $p_r$ is unchanged, and we obtain:
    $$\operatorname*{argmin}_{i=1}^{n=4m-3} \; C_p(\Mcal_i) 
    \;\; = \;\;
    \operatorname*{argmin}_{i=1}^m \; h_i(p).$$ 
    \noindent
    Since consecutive $p_r$ have different unique optimal matchings, the optimal matching changes at least $k-1$ times between points $p_1, \ldots, p_k$; then it changes at least $k-1$ times on $[1/3, 2/3]$.
\end{proof}

\subsection{Proof of the lower bound}

\begin{lemma}\label{lem:lower-bound}
For any $n\geq2$, there exists\footnote{Our proof is designed for $n$ such that $n+3$ is divisible by 4. For other $n$ the result can be achieved via monotonicity; see Lemma \ref{lem:V-0-1-increasing} } a configuration of points and constant $c_1$, such that
$c_1 \cdot n\alpha(n) \leq V_n^*$.
\end{lemma}
\begin{proof}
This fold for $n = 2$ from Lemma \ref{lem:V-0-1-increasing}. For $n\geq 3$, first, use construction from Lemma \ref{lem:analytic-points-in-area} to construct $m = \lfloor\frac{n+3}{4}\rfloor$ segments $s_1, \ldots, s_m$ with $k$ intervals of consecutively different minimisers centred at $p_1 < \cdots < p_k$ and fix $\delta>0$ from that lemma. From that lemma, $k \geq cm\alpha(m) \geq c_1n\alpha(n)$ as $m \approx n/4$ and $\alpha(n)$ is sublinear.
Given this, approximate segments up to $\delta/2$ with functions of type $h_i(p)$ using Lemma \ref{lem:exponential-approximation}. Since minimisers on $(p_r-\delta, p_r + \delta)$ are at least $2\delta$ lower than other segments, after this approximation, at points $p_r$ the minimisers of approximated segments $h_i$ would remain the same, and will be at least $\delta$ apart from other $h$ functions.

This guarantees the conditions for Lemma \ref{lem:points-as-h-i} and shows that the points can be arranged so that the optimal matching switches $k-1$ times. With $k \geq c_1n\alpha(n)$, this proves the bound.
\end{proof}

\section*{Acknowledgements}
The author thanks Andrea Ottolini for the guidance and encouragement throughout this work, and Nicola Guglielmi, Stefano Olla, Filippo Santambrogio, and Aleh Tsyvinski for helpful discussions. \\ 
The author thanks Ypatia Laboratory for Mathematical Sciences (LYSM) for funding a research visit to the Institut Camille Jordan, Lyon, and Center for Algorithms, Data, and Market Design (CADMY) for funding a research visit to Yale University.

{\small
\bibliographystyle{plain}
\bibliography{references}
}

\end{document}